\documentclass[10pt,twoside,english]{amsart}
\advance\oddsidemargin by -1.0cm
\advance\evensidemargin by -1.0cm
\advance\topmargin by -1.0cm 
\usepackage{amssymb}
\usepackage{babel}
\usepackage{amstext}
\usepackage{amscd}   
\usepackage{epsfig}  
\usepackage{rotating}
\theoremstyle{plain}
\newtheorem{cor}{Corollary}[section]
\newtheorem{lem}{Lemma}[section]
\newtheorem{thm}{Theorem}[section]
\newtheorem{prop}{Proposition}[section]
\theoremstyle{definition}

\newcommand{\bdm}{\begin{displaymath}}
\newcommand{\edm}{\end{displaymath}}
\newcommand{\be}{\begin{equation}}
\newcommand{\ee}{\end{equation}}
\newcommand{\ba}[1]{\begin{array}{#1}}
\newcommand{\ea}{\end{array}}

\newcommand{\btab}{\begin{tabular}}
\newcommand{\etab}{\end{tabular}}

\newcommand{\ox}{\otimes}

\newcommand{\lra}{\longrightarrow}

\newcommand{\lan}{\left\langle}
\newcommand{\ran}{\right\rangle}

\renewcommand{\div}{\ensuremath{\mathrm{div}}}

\newcommand{\R}{\ensuremath{\mathbb{R}}}

\newcommand{\Ric}{\ensuremath{\mathrm{Ric}}}

\begin{document}
%
% Die lange Liste der Formeln 1-7 wird als Liste, nicht als equations
% gesetzt !!! --> Trick, um ab da richtige Numerierung zu haben:
% 
\setcounter{equation}{7}
%
%------ draw title page -----
%
\thispagestyle{empty}
%
%------------------------------
\date{\today}
%------------------------------
\title{Eigenvalue estimates for the Dirac operator depending on the 
Weyl tensor}
%------------------------------
%
% author and address
%
%------------------------------
%
\author{Thomas Friedrich}
\author{Klaus-Dieter Kirchberg}
%-------------------------------
\address{ 
{\normalfont\ttfamily friedric@mathematik.hu-berlin.de}\newline
{\normalfont\ttfamily kirchber@mathematik.hu-berlin.de}\newline
Institut f\"ur Reine Mathematik \newline
Humboldt-Universit\"at zu Berlin\newline
Sitz: Rudower Chaussee 25\newline
D-10099 Berlin\\
Germany}
%------------------------------------
\thanks{This work was supported by the SFB 288 "Differential geometry
and quantum physics" of the Deutsche Forschungsgemeinschaft.}
%------------------------------------
%------------------------------------
\begin{abstract}
%--------------
We prove new lower bounds for the first eigenvalue of the Dirac operator
on compact manifolds whose Weyl tensor or curvature tensor, respectively, is
divergence free. In the special case of Einstein manifolds, we obtain
estimates depending on the Weyl tensor.
\end{abstract}
%------------
\maketitle
%----------------
\tableofcontents
%----------------
\pagestyle{headings}
%
%
%-------------- body of the document ------------------------------------------
%
%---------------------------------------------------------------------------- 
\section{Introduction}\noindent
%----------------------------------------------------------------------------
If $M^n$ is a compact Riemannian spin manifold with 
positive scalar curvature $R$, then each eigenvalue $\lambda$ of the Dirac 
operator $D$ satisfies the inequality
\bdm
\lambda^2 \ge \frac{n R_0}{4(n-1)} \ ,
\edm
where $R_0$ is the minimum of $R$ on $M^n$. The estimate is sharp
in the sense that there exist manifolds for which the lower bound is an 
eigenvalue $\lambda_1$ of $D$. In this case $M^n$ must be an Einstein space (see \cite{FR1}). A generalization of this inequality was proved in the paper 
\cite{Hi}, in which a conformal lower bound for the spectrum of the Dirac 
operator appeared. Moreover, for special Riemannian manifolds, better 
estimates 
for the eigenvalues of the Dirac operator are known, see \cite{K1}, \cite{KSW}.
In the paper \cite{FK} we proved an estimate for the eigenvalues of the
Dirac operator depending on the Ricci tensor in case that the curvature tensor
is harmonic. In this note we will prove an estimate of the Dirac spectrum depending on the scalar curvature and on the Weyl tensor for manifolds with divergence free Weyl tensor. In particular, we prove that a compact, conformally Ricci-flat manifold with certain nontrivial 
conformal invariant  $\nu_0$ does not admit any harmonic spinors. A second application of our result refers to symmetric spaces of compact type. In this case our inequality may be simplified and depends mainly on the scalar curvature and on the length of the Weyl tensor. Under the assumption that the curvature tensor is harmonic we prove in section $4$ an estimate depending on the Ricci tensor as well as on the Weyl tensor. \\

%---------------------------------------------------------------------------- 
\section{Curvature endomorphisms of the spinor bundle}\noindent
%----------------------------------------------------------------------------
%
Let $M^n$ be a Riemannian spin manifold of dimension $n \ge 4$ with Riemannian
metric $g$ and spinor bundle $S$. By $\nabla$ we denote the covariant
derivative induced by $g$ on the tangent bundle $TM^n$ as well as the 
corresponding derivative in the spinor bundle $S$. For any vector fields 
$X,Y$ on $M^n$, we use the notation
\bdm
\nabla_{X,Y}\ :=\ \nabla_X\nabla_Y - \nabla_Y\nabla_X
\edm
for the tensorial derivatives of second order. By $K$ and $C$ we denote
the Riemannian curvature tensor and the curvature tensor of $S$,
respectively. Then, for any vector fields $X,Y,Z$ and any
spinor field $\psi$, we have 
\bdm
K(X,Y)(Z)\ =\ \nabla_{X,Y}Z -\nabla_{Y,X}Z,\quad
C(X,Y) \psi \ =\ \nabla_{X,Y} \psi - \nabla_{Y,X} \psi \,.
\edm
Given a local frame of vector fields $(X_1,\ldots,X_n)$, we denote by
$(X^1,\ldots,X^n)$ the associated frame defined by $X^k:=g^{kl}X_l$,
where $(g^{kl})$ is the inverse of the matrix $(g_{kl})$ 
with $g_{kl} = g(X_k,X_l)$. 
For the reader's convenience, we summarize some well known identities:
\begin{enumerate}
\item\label{gl1} $C(X,Y)\,=\, \frac{1}{4}X_k\cdot K(X,Y)(X^k)\,= - 
\frac{1}{4} K(X,Y)(X^k) \cdot X_k$,
\vspace{1mm}
\item\label{gl2} $C(X,Y) \cdot Z - Z \cdot C(X,Y)=K(X,Y)(Z)$,
\vspace{1mm}
\item\label{gl3} $X_k \cdot C(X^k,Y) \ = \ \frac{1}{2} {\rm Ric} (Y) \ = 
\ C(X^k,Y) \cdot X_k$,
\vspace{1mm} 
\item\label{gl4} $X_k \cdot K(X,X^k)(Y) \ = \ 2C(X,Y)+g (\Ric (X),Y)$,
\vspace{1mm}
\item\label{gl5} $X^k \cdot \nabla_{X,X_k} \psi \ = \ \nabla_X D \psi$,
\vspace{1mm}
\item\label{gl6} $X^k \cdot \nabla_{X_k ,X} \psi \ = \ \nabla_X D \psi+
\frac{1}{2} \Ric (X) \cdot \psi$,
\vspace{1mm}
\item\label{gl7} $X^k \cdot \Ric (X_k) \ = \ \Ric (X_k) \cdot X^k \ =
\ -R$,
\end{enumerate}

where $\Ric$ is the Ricci tensor, $R$ the scalar curvature and $D$ the
Dirac operator locally defined by $\Ric (X) := K(X,X_k)(X^k), R=g
(\Ric(X_k), X^k)$ and $D\psi = X^k \cdot \nabla_{X_k} \psi$, respectively.
For any vector fields $X$ and $Y$, we consider the endormorphism $E(X,Y)$ 
of $S$ or $\Gamma (S)$, respectively, locally given by
\bdm
E(X,Y) \ := \ - \, C(X_k,Y) \cdot C(X^k,X) \, . 
\edm
Since $C(X,Y)$ is anti-selfadjoint with respect to the Hermitian scalar
product $\langle \cdot , \cdot \rangle$ on $S$, the endomorphism $E(X,Y)$
has a similar property:
\begin{equation}\label{gl8}
C(X,Y)^* \ = \ - \, C(X,Y), \quad E(X,Y)^* \ = \ E(Y,X) \, . 
\end{equation}
The endomorphism $F$ of $S$ defined as the contraction of $E$,
\bdm
F \  := \ E(X_k , X^k) \ = \ - \, C(X_k,X_k)C(X^k,X^l) \ ,
\edm
is selfadjoint and nonnegative, i.e., 
\begin{equation}\label{gl10}
F^* \ = \ F  , \quad \langle F \psi , \psi \rangle \ \ge \ 0 \, . 
\end{equation}
We denote by $W$ the Weyl tensor of the Riemannian manifold and we introduce 
the following endomorphisms acting in the spinor bundle:
\begin{eqnarray*}
B(X,Y) &:=& \frac{1}{4} X_k \cdot W(X,Y)(X^k) \ , \quad
G(X,Y) :=  - \, B(X_k,Y)B(X^k,X) \ , \\[0.5em]
H&:=& G(X_k,X^k) \ = \ - \, B(X_k, X_l)B(X^k,X^l) \, .
\end{eqnarray*}
We collect some properties of these endomorphisms:
\begin{equation}\label{gl11}
B(X,Y)^* \ = \ -B(X,Y) , \quad G(X,Y)^* \ = \ G(Y,X)
\end{equation}
\begin{equation}\label{gl12}
X_k \cdot B(X^k,Y) \ = \ 0 \ = \ B(X^k ,Y) \cdot X_k \ , 
\end{equation}
\begin{equation}\label{gl13}
H^* \ = \ H  , \quad \langle H \psi , \psi \rangle \ \ge \ 0 \ , 
\end{equation}
\begin{equation}\label{gl14}
B(X,Y) \cdot Z - Z\cdot B(X,Y) \ = \ W(X,Y)(Z) \ , 
\end{equation}
\begin{equation}\label{gl16}
F \ = \ H+ \frac{1}{2(n-2)} \Big|\Ric - \frac{R}{n} \Big|^2 + \frac{R^2}{4n(n-1)}
\ . 
\end{equation}
Let $\nu(x)$ denote the smallest eigenvalue of $H$ at the point $x \in M^n$. 
Its infimum
\bdm
\nu_0 \ := \ \inf \{ \nu(x) \ : x \in M^n \} \ \ge \ 0
\edm
will occur in our estimates of the eigenvalues of the Dirac operator $D$.
If $\psi$ is a parallel spinor $(\nabla \psi =0)$, then, for all vector
fields $X$ and $Y$, it follows that $C(X,Y) \cdot \psi = 0$ and, hence, $F \psi =0$.
Thus, the relation (\ref{gl16}) shows that the Ricci tensor as well as
the function $\nu \ge 0$ are obstructions for the existence of parallel
spinors. Moreover, the Schr\"odinger-Lichnerowicz formula
\begin{equation}\label{gl17}
\nabla^* \nabla \ = \ D^2 - \frac{R}{4}
\end{equation}
implies that, on compact manifolds with vanishing scalar curvature, each 
harmonic spinor is parallel. Hence, as $M^n$ is compact, Ricci flat and
$\nu_0 >0$, there are no harmonic spinors.
%
%---------------------------------------------------------------
\section{Estimate for manifolds with divergence free Weyl tensor}
%---------------------------------------------------------------
In this section we assume that the Weyl tensor $W$ of $M^n$ is
divergence free, i.\,e., for all vector fields
$Y$ and all local frames  $(X_1,\ldots,X_n)$ $W$ satisfies the condition
\begin{equation}\label{gl19}
(\nabla_{X_k}W)(X^k,Y)\ =\ 0\,.
\end{equation}
By definition of $B$, this implies  
\begin{equation}\label{gl20}
(\nabla_{X_k}B)(X^k,Y) = 0 \ . 
\end{equation}
For any real parameter $t\in\R$, we consider the differential operator 
$P^t:\ \Gamma(S)\lra \Gamma(TM^n\ox S)$ locally defined by
$P^t\psi:=X^k\ox P^t_{X_k}\psi$ and
\begin{equation}\label{gl21}
P^t_{X}\psi\ :=\ \nabla_X\psi+\frac{1}{n}X\cdot D\psi - tB(X,X^k)\cdot
\nabla_{X_k}\psi\,.
\end{equation}
Using the twistor operator $\mathcal{D}: \Gamma (S) \lra \Gamma (TM^n \ox S)$
given by $\mathcal{D} \psi := X^k \ox ( \nabla_{X_k}\psi + \frac{1}{n} X_k \cdot D \psi)$ this may be rewritten as
\begin{equation}\label{gl211}
P^t\psi\ :=\ \mathcal{D} \psi - t X^k \ox B(X_k,X^l)\cdot
\nabla_{X_l}\psi\,.
\end{equation}
The image of $\mathcal{D}$ is contained in the kernel of the Clifford multiplication, i.\,e., it holds that
\begin{equation}\label{gl22}
X^k \cdot \mathcal{D}_{X_k} \psi \ = \ 0 \, . 
\end{equation}
Thus, by (\ref{gl12}) and (\ref{gl22}), it follows that
\begin{equation}
X^k \cdot P^t_{X_k} \psi \ = \ 0 \, . 
\end{equation}
\begin{lem}\label{Pt2}
%---------------------
Suppose that the Weyl tensor $W$ is divergence free. Then, any spinor
field $\psi$ satisfies the equation
\begin{equation}\label{eq-Pt2}
|P^t\psi|^2\ =\ |\nabla\psi|^2 -\frac{1}{n}|D\psi|^2- t \lan H \psi, 
\psi \ran + t^2 \lan G(X^k,X^l) \nabla_{X_k}\psi , \nabla_{X_l} \ran 
- 2t\, \div \lan B \psi, \nabla \psi \ran  \, ,
\end{equation}
where $\lan B\psi,\nabla\psi\ran$ is the vector field locally
defined by 
$\lan B\psi,\nabla\psi\ran:=\lan B(X^i,X^k)\psi,\nabla_{X_k}\psi\ran X_i$.
\end{lem}
\begin{proof}
%------------
Using the formulas (\ref{gl11}), (\ref{gl12}) and (\ref{gl21}) we calculate
\begin{eqnarray*}
|P^t \psi |^2 &=& \langle P^t_{X_i} \psi , P^t_{X^i} \psi \rangle \\[0.5em]
&=& |D \psi |^2 - 2t \langle \nabla_{X_i} \psi , B (X^i, X^k)
\nabla_{X_k} \psi \rangle 
+t^2 \langle B(X_i , X^k) \nabla_{X_k} \psi , B(X^i, X^l)
\nabla_{X_l} \psi \rangle \ . 
\end{eqnarray*}
Thus, we obtain
\bdm
(*)\quad \quad
|P^t \psi |^2 \ = \ | \nabla \psi |^2 - \frac{1}{n} |D \psi|^2- 2t \langle
\nabla_{X_i} \psi , B(X^i, X^k) \nabla_{X_k} \psi \rangle 
+ t^2 \langle G(X^k, X^l) \nabla_{X_k} \psi , \nabla_{X_l} \psi
\rangle 	 \ . 
\edm 
Let $x \in M^n$ be any point and let $(X_1 , \ldots , X_n)$ be an
orthonormal frame in a neighbourhood of $x$ with $(\nabla X_k)_x =0$.
Then, we have at the point $x$
\begin{eqnarray*}
\langle \nabla_{X_i} \psi , B(X^i , X^k) \nabla_{X_k} \psi \rangle &=&
X_i (\langle \psi , B(X^i, X^k)\nabla_{X_k} \psi \rangle ) - \langle \psi , (\nabla_{X_i} B)(X^i, X^k) \nabla_{X_k} \psi \rangle \\[0.5em]
&& - \langle \psi , B(X^i, X^k)\nabla_{X_i} \nabla_{X_k} \psi \rangle \\[0.5em]
&\stackrel{(17)}{=}& \div \langle B \psi , \nabla \psi \rangle - \frac{1}{2} 
\langle \psi, B(X^i, X^k)C(X_i , X_k) \psi \rangle \\[0.5em]
&\stackrel{(11)}{=}& \div \langle B \psi , \nabla \psi \rangle - 
\frac{1}{2} \langle \psi, B(X^i, X^k)B(X_i , X_k) \psi \rangle \\[0.5em]
&=& \div \langle B \psi , \nabla \psi \rangle + \frac{1}{2} \langle \psi, 
H \psi \rangle \ . 
\end{eqnarray*}
Inserting this into $(*)$ we obtain (\ref{eq-Pt2}).
\end{proof}
Let us introduce the number $\mu_0$ measuring the maximum of the norm of the
Weyl tensor,
\bdm
\mu_0^2 \ = \ \max \Big( \frac{1}{16}||W_{XYij} X^i \cdot X^j||^2 : 
 \ \lan X,Y \ran = 0 , \ |X| = |Y| = 1 \Big)
\edm
Then, for any point $x \in M^n$ and any orthonormal basis $(X_1, \ldots , X_n)$
of $T_x M^n$, we have
\bdm
\hspace{-3cm}|\langle G(X^k, X^l)\nabla_{X_k} \psi , \nabla_{X_l} \psi \rangle | \ = \ | \langle B(X_i , X^k) \nabla_{X_k} \psi , B(X^i , X^l) \nabla_{X_l}
\psi \rangle |
\edm
\bdm
\le \ \sum\limits_{i,k.l} |\langle B(X_i , X_k) \nabla_{X_k} \psi,
B(X_i , X_l) \nabla_{X_l} \psi \rangle |
\ \le \ \sum\limits_{i,k,l} | B(X_i , X_k) \nabla_{X_k} \psi | \cdot |
B(X_i , X_l) \nabla_{X_l} \psi|
\edm
\bdm
\hspace{-1.7cm}\le \ \sum\limits_{i,k,l} \| B(X_i , X_k) \| \, \| B(X_i, X_l) 
\| \, | \nabla_{X_k} \psi | \, |\nabla_{X_l} \psi |
\ \le \  n \mu^2_0 \Big( \sum\limits_{k,l} |\nabla_{X_k} \psi | \, |\nabla_{X_l}
\psi | \Big) 
\edm
\bdm\ 
\hspace{-3.9cm} = \ n \mu^2_0 \Big( \sum\limits_k | \nabla_{X_k} \psi | \Big)^2 \
\le \  n^2 \mu^2_0 \Big( \sum\limits_k | \nabla_{X_k} \psi |^2 \Big) \ = \
n^2 \mu^2_0 | \nabla \psi |^2 \ . 
\edm
Thus, we obtain the estimate
\begin{equation}\label{gl25}
| \langle G(X^k , X^l)\nabla_{X_k} \psi , \nabla_{X_l} \psi \rangle | 
\le n^2 \mu_0^2 |\nabla \psi |^2 \ . 
\end{equation}
Furthermore, let us denote the minimum of $R$ on $M^n$ by $R_0$. Then,
applying (\ref{gl17}), (\ref{eq-Pt2}), (\ref{gl25}) to an eigenspinor $\psi$ 
of the Dirac operator $(D \psi = \lambda \psi)$ we have
\begin{eqnarray*}
0 &\le & \int_{M^n} |P^t \psi |^2 \ = \ \int_{M^n} \Big(
\frac{n-1}{n} \lambda^2 | \psi |^2 - \frac{R}{4} |\psi|^2 - 
t \langle H \psi, \psi \rangle + t^2 \langle G(X^k, X^l) \nabla_{X_k} \psi ,
\nabla_{X_l} \psi \rangle \Big)\\[0.5em]
&\le& \Big( \frac{n-1}{n} \lambda^2 - \frac{R_0}{4} - \nu_0 t + n^2 \mu_0^2 
\big( \lambda^2 - \frac{R_0}{4} \big) t^2 \Big) \cdot 
\int_{M^n} | \psi |^2 \ . 
\end{eqnarray*}
This implies the inequality
\bdm
\frac{n-1}{n} \lambda^2 - \frac{R_0}{4} - \nu_0 t + n^2 \mu_0^2 (\lambda^2
- \frac{R_0}{4} )\, t^2 \ge 0 \ ,
\edm
which is equivalent to
\begin{equation}\label{gl26}
\lambda^2 \ge \frac{n}{4} \cdot \frac{R_0 + 4\nu_0 t +n^2 R_0 \mu_0^2 
t^2}{(n-1)+ n^3 \mu_0^2  t^2} \quad \quad (\forall t \in \R) \ . 
\end{equation}
By computing the maximum of the right side with respect to the parameter $t$ 
we obtain the main theorem of this section.
\begin{thm}\label{main-thm}
%--------------------------
Let $M^n$ be a compact Riemannian spin manifold with divergence 
free Weyl tensor. Then, for any eigenvalue $\lambda$ of the
Dirac operator, we have the inequality
\begin{equation}\label{gl27}
\lambda^2 \ge \frac{n R_0}{4(n-1)} + \frac{2 \nu_0^2}{n \mu_0^2
\Big( R_0 + \sqrt{R_0^2 + \frac{n-1}{n} (\frac{4\nu_0}{\mu_0})^2}\Big)}
\, . 
\end{equation}
Moreover, for $R_0 <0$, this lower bound is positive if the condition
$\nu_0 > \frac{n}{2} | R_0| \, \mu_0$ is satisfied. 
\end{thm}
The divergence of the Weyl tensor is given by the well known identity
\begin{equation}\label{gl29}
(\nabla_{X_k} W)(X,Y)(X^k) \ = \ (n-3)\Big(( \nabla_X T)(Y)-(
\nabla_Y T)(X) \Big) \ , 
\end{equation}
where the tensor $T$ is defined by
\bdm
T(X) \ := \ \frac{1}{n-2} \Big( \frac{R}{2(n-1)} \, X - {\rm Ric} (X) \Big)
\ . 
\edm
In particular, any Einstein manifold has a divergence free Weyl tensor. 
\begin{cor}\label{cor-3}
%-----------------------
The inequality  ($\ref{gl27}$) holds for any compact Einstein spin manifold. 
\end{cor}
\begin{cor}\label{cor-4}
%-----------------------
Let $M^n$ be a compact Riemannian spin manifold with divergence
free Weyl tensor and vanishing scalar curvature. Then we have the estimate
\begin{equation}\label{gl30}
\lambda^2 \ge \frac{\nu_0}{2 \mu_0 \sqrt{n(n-1)}} \ . 
\end{equation}
\end{cor}
Consider a compact Riemannian spin manifold $(M^{n}, g^*)$ and suppose that 
there exists a conformally equivalent Ricci-flat metric $g = e^f \cdot g^*$. 
Since the Weyl tensor is a conformal invariant, the condition $\nu_0 > 0$ is 
conformally invariant, too. The same is true for the dimension of the space 
of harmonic spinors.
\begin{cor}\label{cor-}
A compact, conformally Ricci-flat spin manifold with $\nu_0 > 0$ does not admit harmonic spinors.
\end{cor}
In case of an even-dimensional manifold the spinor bundle splits into the bundle of positive and negative spinors, respectively. We can introduce two smallest eigenvalues $\nu_0^{\pm}$ and Corollary 3.3 holds in any of these two bundles. \\

The lower bound for the eigenvalues of the Dirac operator proved in Theorem 3.1
depends on the minimum of the scalar curvature, on the maximum $\mu_0^2$ 
and on the smallest eigenvalue of the non-negative endomorphism
\bdm
H \ := \ - \frac{1}{16} \sum_{k,l} \sum_{\alpha, \beta,  \gamma, \delta}
W_{kl\alpha\beta}W_{kl\gamma\delta} X^{\alpha} \cdot X^{\beta} \cdot X^{\gamma}
\cdot X^{\delta} 
\edm
acting on the spinor bundle. Using the grading of the Clifford algebra we
decompose the endomorphism $H := H_0 + H_2 + H_4$ into three parts, where
$H_0$ is a scalar and $H_2,H_4$ are elements of the Clifford algebra of
degree two and four, respectively. It is easy to compute $H_0$,
\bdm
H_0 \ = \ \frac{1}{8} |W|^2 \ . 
\edm
Since $H$ and $H_4$ are hermitean and $H_2$ is anti-hermitean, we conclude
that $H_2 = 0$. Consequently, we obtain the formula
\bdm
\nu_0 \ = \ \min \Big( \frac{1}{8}|W|^2 + \lan H_4 \cdot \psi, \psi \ran \, :
\ |\psi| = 1 \Big) \ ,
\edm
where $H_4$ is given by 
\bdm
H_4 \ = \ - \frac{1}{2} \sum_{k,l} \sum_{\alpha< \beta< \gamma<\delta}
\Big(W_{kl\alpha\beta}W_{kl\gamma\delta}- W_{kl\alpha\gamma}W_{kl\beta\delta}
+ W_{kl\alpha\delta}W_{kl\beta\gamma} \Big) X^{\alpha} \cdot X^{\beta} 
\cdot X^{\gamma} \cdot X^{\delta} \ .
\edm
In the four-dimensional case the endomorphism $X^1 \cdot X^2 \cdot X^3 
\cdot X^4$ acts on the two parts of the spinor bundle $S = S^+ \oplus S^-$ by 
multiplication by $\pm 1$ and we obtain the simple formula
\bdm
\nu_0 \ = \ \min \Big( \frac{1}{4}|W + *W|^2 ,\, \frac{1}{4}|W - *W|^2 \Big) \ .
\edm
For example, consider the square of the Dirac operator acting on the bundle
$S^+$ of positive spinors over a $4$-dimensional K\"ahler-Einstein manifold.
The positive part $W^+$ of the Weyl tensor acting on $\Lambda^2_+$ has the
diagonal form $W^+ = \mbox{diag}( - R/6, \, - R/6, \, R/3)$ 
(see \cite{Besse}) and, consequently, we compute
\bdm
 \nu_0^+ \ = \ \frac{1}{6} R^2, \quad \mu_0^2 \ = \ \frac{1}{8 \cdot 9} R^2 \ .
\edm
The estimate of Theorem 3.1 yields the inequality $\lambda^2 \ge R/2$ and this is precisely the lower bound for the eigenvalues of the Dirac operator
on any $4$-dimensional K\"ahler manifold (see  \cite{K1} and \cite{FR3}).\\

On the other hand, if the Weyl tensor of the manifold $M^n \, (n \ge 5)$ satisfies the relation
\bdm
 \sum_{k,l} \Big(W_{kl\alpha\beta}W_{kl\gamma\delta} + W_{kl\alpha\gamma}
W_{kl\delta\beta} + W_{kl\alpha\delta}W_{kl\beta\gamma} \Big) \ = \ 0 \, ,
\edm
then $H_4$ vanishes. This situation occurs, for example, if $M^n$ is an 
irreducible symmetric space of compact type. The proof is an easy computation 
using the well-known formulas for the curvature tensor of a symmetric space,
which is why we shall only sketch it shortly. First we remark that the following relation holds on any Einstein space:
\bdm
\sum_{k,l}W_{kl\alpha\beta}W_{kl\gamma\delta} \ = \ \sum_{k,l}R_{kl\alpha\beta}R_{kl\gamma\delta} + \frac{4R}{n(n-1)} R_{\alpha\beta\gamma\delta} + \frac{2R^2}{n^2(n-1)^2} \Big(\delta_{\alpha\gamma}\delta_{\beta\delta} - \delta_{\alpha\delta}\delta_{\beta\gamma}\Big) \, .
\edm
In case of a symmetric space $M^n=G/K$ we obtain the formula
\bdm
 \sum_{k,l}R_{kl\alpha\beta}R_{kl\gamma\delta} \ = \ \big(\Omega^G - \Omega^K \big) R_{\alpha\beta\gamma\delta} \ ,
\edm 
where $\Omega^G$ and $\Omega^K$ are the Casimir operators acting on the Lie
algebras of the group $G$ and $K$, respectively. The result follows now from 
the first Bianchi identity of the curvature tensor.
\begin{prop}\label{prop-1}
%------------------------
Let $M^n$ be an irreducible symmetric space of compact type. Then
\bdm
\nu_0 \ = \ \frac{1}{8} |W|^2 \ .
\edm 
\end{prop}
It is well-known that a compact, symmetric space with $\lambda^2 = 
\frac{nR}{4(n-1)} > 0$ is a sphere. This result follows immediately from our
inequality. Indeed, in this case $M^n$ is an irreducible Einstein manifold
(see \cite{FR1}, \cite{FR2}) and, furthermore, we have $\nu_0 = 0$. Since the 
manifold is symmetric, we conclude that the Weyl tensor vanishes, i.e., $M^n$ 
is a space of constant curvature.
\section{Estimates in case of divergence free curvature tensor}

In this section we assume that the Riemannian curvature tensor $K$
is divergence free,i.e., locally we have the equality
\begin{equation}\label{gl39}
(\nabla_{X_k} K)(X^k , Y) \ = \ 0 
\end{equation}
for each vector field $Y$. The Bianchi identity implies the general
relation
\begin{equation}\label{gl40}
(\nabla_{X_k} K)(X,Y)(X^k) \ = \ (\nabla_Y \Ric)(X) - (\nabla_X \Ric)(Y)
\ . 
\end{equation}

Thus, (\ref{gl39}) is equivalent to
\begin{equation}\label{gl41}
(\nabla_X \Ric)(Y) \ = \ (\nabla_Y \Ric)(X) \ . 
\end{equation}

In particular, the scalar curvature $R$ is constant. Moreover, (\ref{gl39})
implies
\begin{equation}\label{gl42}
(\nabla_{X_k} C)(X^k , Y) \ = \ 0 \ . 
\end{equation}

For $t \in \R$, we now consider the operator
\bdm
Q^t : \Gamma (S) \rightarrow \Gamma (TM^n \otimes S)
\edm

defined by $Q^t \psi := X^k \otimes Q^t_{X_k} \psi$ and
\bdm
Q^t_X \psi := {\mathcal D}_X \psi - t \cdot C(X,X^k) \cdot \nabla_{X_k} \psi
\ . 
\edm

A straightforward calculation yields
\begin{eqnarray}\label{gl43}
|Q^t \psi |^2  &=&  | {\mathcal D} \psi |^2 + 2 t \langle C(X^k, X^l)
\nabla_{X_k} \psi , \nabla_{X_l} \psi \rangle 
 + \frac{t}{n} {\rm Re} \langle  D \psi, \Ric (X^k)
\nabla_{X_k} \psi \rangle \\[1mm]
\notag &+& t^2 \langle E(X^k, X^l) \nabla_{X_k} \psi , 
\nabla_{X_l} \psi \rangle \ . 
\end{eqnarray}

Furthermore, by Lemma 1.2 and Lemma 1.4 in \cite{FK} and our assumption, 
there are the equations
\begin{eqnarray}\label{gl44}
\langle C(X^k, X^l) \nabla_{X_k} \psi , \nabla_{X_l} \psi \rangle
&=& \div \langle C \psi , \nabla \psi \rangle - \frac{1}{2} \langle
\psi , F \psi \rangle \ , \\[0.5em]
\label{gl45}
{\rm Re} \langle \psi, \Ric (X^k) \nabla_{X_k} D \psi \rangle &=& | \nabla
D \psi |^2 - |(D^2 - \frac{R}{4} )\psi |^2 - \frac{R}{4} | \nabla \psi |^2
\\[0.5em]
&& + \frac{1}{4} | \Ric |^2 | \psi |^2 + \langle \nabla_{\Ric (X_k)} \psi,
\nabla_{X^k} \psi \rangle - \div (X_{\psi} ) \ , \nonumber
\end{eqnarray}

where $X_{\psi}$ is a vector field depending on $\psi$. From (\ref{gl16}),
(\ref{gl43}), (\ref{gl44}) and (\ref{gl45}) we obtain the basic Weitzenb\"ock formula of this section.

\begin{lem}\label{lem4-1}
Let $M^n$ be a Riemannian spin manifold with divergence free curvature tensor
and let $\lambda$ be any eigenvalue of the Dirac operator.
Then, for any corresponding eigenspinor $\psi$ and for all $t \in \R$, we
have the equation
\bdm \label{gl46}
| Q^t \psi |^2 \ = \ | \nabla \psi |^2 - \frac{\lambda^2}{n} | \psi |^2
+ \frac{t}{n} \langle \nabla_{\Ric (X_k)} \psi , \nabla_{X^k} \psi
\rangle 
+ \frac{t}{n} (\lambda^2 - \frac{R}{4})\Big(| \nabla \psi |^2 - (\lambda^2 - 
\frac{R}{4})|\psi |^2\Big) - t \langle \psi , H \psi \rangle
\edm
\bdm
- \frac{t}{4n}  \Big( \frac{n+2}{n-2} \Big| \Ric - \frac{R}{n} \Big|^2 +
\frac{R^2}{n(n-1)} \Big) | \psi |^2 + \div \Big(2t \langle C \psi , 
\nabla \psi \rangle - \frac{t}{n} X_{\psi} \Big)
+ t^2 \langle E(X^k , X^l) \nabla_{X_k} \psi , \nabla_{X_l} \psi \rangle
 \ .
\edm
\end{lem}

We consider the curvature tensor $K$ as an endomorphisms of $\Lambda^2 M^n$ by
the usual definition
\bdm
K(u^i \wedge u^j) \ := \ \frac{1}{2} \, g(K(X^i, X^j)(X_k), X_l) u^k \wedge u^l
\ , 
\edm

where $(u^1 , \ldots , u^n)$ denotes the coframe dual to the local 
frame $(X_1 , \ldots , X_n)$. Then $K$ is
selfadjoint  with respect to the scalar product on $\Lambda^2 M^n$ induced by 
the Riemannian metric $g$. Let
$M^n$ be compact and denote by $\sigma$ the maximum of the absolut values of 
the eigenvalues of $K$ on $\Lambda^2 M^n$. Then we estimate the operator norm
of the endomorphism $C(X_i,X_j)$ acting on the spinor bundle
\begin{eqnarray*}
\| C(X_i , X_j) \| &\le & \frac{1}{2} \, \sum\limits_{k<l} | g(K(X_i , X_j)(X_k), X_l)| \cdot \| X_k \cdot X_l \| \\[0.5em]
&=& \frac{1}{2} \, \sum\limits_{k<l} | g (K(X_i , X_j)(X_k), X_l)| \ \le \ 
\frac{1}{2} \Big( \begin{array}{c}n\\2 \end{array} \Big) \sigma \ . 
\end{eqnarray*}

\newcommand{\kleh}{\Big( \begin{array}{c}n\\2 \end{array} \Big)}

Using this upper bound we obtain
\bdm
\Big| \langle E(X^k , X^l) \nabla_{X_k} \psi , \nabla_{X_l} \psi \rangle \Big| 
 \ \le \ \sum\limits_{i,k,l} \| C(X_i , X_k)\| \cdot \| C(X_i , X_l)\| \cdot 
| \nabla_{X_k} \psi | \cdot | \nabla_{X_l} \psi | 
\edm
\bdm
\le \frac{n}{4} \kleh^2 \sigma^2 \Big( \sum\limits_{k,l} | \nabla_{X_k} \psi |
| \nabla_{X_l} \psi | \Big) 
\ \le \  \frac{n^2}{4} \kleh^2 \sigma^2 \Big( \sum\limits_k | \nabla_{X_k}
\psi |^2 \Big) \ = \ \Big( \frac{n}{2} \kleh \sigma \Big)^2 | \nabla \psi
|^2  
\edm

and, consequently, the inequality
\begin{equation}\label{gl48}
\Big| \langle E(X^k, X^l) \nabla_{X_k} \psi , \nabla_{X_l} \psi \rangle
\Big| \le \Big( \frac{n}{2} \kleh \sigma \Big)^2 | \nabla \psi |^2 \ . 
\end{equation}

If $\kappa$ denotes the maximum of all eigenvalues of the Ricci tensor on
$M^n$, then we have 
\begin{equation}\label{gl49}
\langle \nabla_{\Ric (X_k)} \psi , \nabla_{X^k} \psi \rangle \le \kappa
| \nabla \psi |^2 \ . 
\end{equation}

Then, by Lemma (\ref{lem4-1}), (\ref{gl48}) 
and (\ref{gl49}), for all $t \ge 0$, we obtain
\begin{eqnarray}\label{gl50}
\Big( \frac{n-1}{n} + \frac{\kappa}{n} t+ \Big( \frac{n}{2} \kleh 
\sigma \Big)^2 t^2 \Big) \lambda^2  \ge 
&+& \frac{1}{4n} \Big( \frac{n+2}{n-2} \Big| \Ric - \frac{R}{n} \Big|^2_0 +
\frac{R^2}{n(n-1)} \Big) t + \nu_0 t \\[1mm]
&+& \frac{R}{4} \Big( 1+ \frac{\kappa}{n} t + \Big( \frac{n}{2} \kleh
\sigma \Big)^2 t^2 \Big)\nonumber \ ,
\end{eqnarray}

where $|\Ric - \frac{R}{n} |_0$ denotes the minimum of the length on $M^n$.
Inserting $\lambda =0$ in this inequality we obtain

\begin{thm}\label{thm-4-1}
Let $M^n$ be a compact Riemannian spin manifold with divergence free
curvature tensor and scalar curvature $R \le 0$ such that the condition
\begin{equation}\label{gl51}
\frac{n+2}{n-2} \Big| \Ric - \frac{R}{n} \Big|^2_0 + \frac{R^2}{n(n-1)} +
4n \nu_0 > |R| \Big( \kappa + n^2 \kleh \sigma \Big) 
\end{equation}
is satisfied. Then there are no harmonic spinors.
\end{thm}
For simplicity, let us introduce the notations
\bdm
a \ := \ 4n (n-1) \nu_0 + \frac{(n-1)(n+2)}{n-2} \Big| \Ric - \frac{R}{n} 
\Big|^2_0 - R (\kappa - \frac{R}{n} ) \, ,
\edm
\bdm
b \ := \ \frac{1}{2} n^2 \kleh^3 \sigma^2 \quad , \quad A_{\pm} \ := \
\sqrt{b^2 R^2 + ab(a+R\kappa)} \pm bR \ . 
\edm
Moreover, using the new parameter $s: = \frac{t}{n-1} \ge 0$, the
inequality (\ref{gl50}) can be written as
\begin{equation}\label{gl52}
\lambda^2 \ge \frac{nR}{4(n-1)} + \frac{s}{4(n-1)} \cdot
\frac{a-bRs}{1+ \kappa s+bs^2}  \ . 
\end{equation}
Calculating the maximum of the right side with respect to $s\ge 0$ we obtain 
our main result.
\begin{thm}\label{thm-4-2}
Let $M^n$ be a non-flat compact Riemannian spin manifold with 
divergence free curvature tensor and let $\lambda$ be any eigenvalue of the Dirac operator. Then we have the estimate
\begin{equation}\label{gl53}
\lambda^2 > \frac{nR}{4(n-1)} + \frac{a}{4(n-1)} \cdot \frac{A_-}{2ab+
\kappa A_+}  \ . 
\end{equation}
\end{thm}
Remark that in case $R>0, a>0$ this lower bound is greater then
$nR/4(n-1)$. If $R \le 0$ the lower bound is positive under the condition
\begin{equation}\label{gl54}
a + R \kappa > (n-1) |R| \Big( \kappa+n^2 \kleh \sigma \Big) \ .
\end{equation}
\begin{proof} 
It remains to show that, for the first eigenvalue $\lambda_1$ of
$D$, in (\ref{gl53}) equality can not occur. Let us assume the counterpart. 
Then any eigenspinor $\psi$ corresponding to $\lambda_1$ satisfies the
equation $Q^{t_0} \psi =0$ with the optimal parameter $t_0 >0$. By (\ref{gl3}) and (\ref{gl22}), $Q^{t_0} \psi =0$ implies\\

$(*) \hfill \Ric (X^k) \nabla_{X_k} \psi = 0 \ . \hfill \mbox{}$\\

Moreover, the limiting case of (\ref{gl53}) implies that in the inequality
for $\| C(X_i, X_j )\|$  we have an equality, i.e., 
\bdm
\| C(X_i, X_j )\| = \frac{1}{2} \kleh \sigma  \ .
\edm
Hence, $M^n$ is a space of constant sectional curvature. In particular,
$M^n$ is Einstein $(\Ric = \frac{R}{n})$ and $(*)$ implies  
$0= R \cdot D \psi = R \cdot \lambda_1 \cdot \psi$. Consequently, the Ricci
tensor vanishes and $M^n$ is flat, a contradiction.
\end{proof}
The curvature tensor of any Einstein manifold of dimension $n \ge 4$ is 
divergence free. In this special case we have $\kappa = R/n$ and the number 
$a$ simplifies to $a =4n (n-1) \nu_0$.\\ 

Finally we remark that we  generalized Corollary 3.2.
\begin{cor}\label{cor-4-1}
Let $M^n$ be a compact Riemannian spin manifold with divergence
free curvature tensor and vanishing scalar curvature such that at least
one of the numbers $\nu_0$ or $|\Ric|_0$, respectively, is not zero.
Then all eigenvalues $\lambda$ of the Dirac operator satisfy the inequality
\bdm
\lambda^2 > \frac{(n+2) |\Ric|_0^2 +4n(n-2) \nu_0}{4(n-2) (\kappa +
\kleh \sigma \sqrt{n(n-1)})} \ . 
\edm
\end{cor}
%
%
%
%-----------------------------------------------------------------------

%


\begin{thebibliography}{1111}
%-----------------------------------------------------------------------
\bibitem{Besse}
A.L.~Besse, \emph{Geometrie riemannienne en dimension $4$}, Cedic-Fernand  
Nathan, Paris 1981.
\bibitem{FR1}
Th.~Friedrich, \emph{Der erste Eigenwert des Dirac-Operators einer kompakten
Riemannschen Mannigfaltigkeit nichtnegativer Skalarkr\"ummung}, Math. Nachr.
97 (1980), 117-146.
\bibitem{FR2}
Th.~Friedrich, \emph{A remark on the first eigenvalue of the Dirac operator
on $4$-dimensional manifolds}, Math. Nachr. 102 (1981), 53-56.
\bibitem{FR3}
Th.~Friedrich, \emph{The classifiaction of $4$-dimensional K\"ahler manifolds
with small eigenvalue of the Dirac operator}, Math. Ann. 295 (1993), 565-574.
\bibitem{FK}
Th.~Friedrich and K.-D.~Kirchberg, \emph{Eigenvalue estimates of the Dirac
operator depending on the Ricci tensor}, to appear.
\bibitem {Hi} 
O. ~Hijazi, \emph{A conformal lower bound for the smallest eigenvalue
of the Dirac operator and Killing spinors}, Commun. Math. Phys. 104 (1986), 
151-162.
\bibitem{K1}
K.-D.~Kirchberg, \emph{The first eigenvalue of the Dirac
operator on K\"ahler manifolds}, J. Geom. Phys. 7 (1990), 449-468.
\bibitem {KSW} 
W.~Kramer, U.~Semmelmann, G.~Weingart, \emph{Eigenvalue estimates
for the Dirac operator on quaternionic K\"ahler manifolds}, Math. Z.
230 (1999), 727-751.


%--------------------
\end{thebibliography}
\end{document}